\input amsppt.sty
\magnification=\magstep1
\hsize=30truecc
\baselineskip=16truept
\vsize=22.2truecm
\nologo
\pageno=1
\topmatter
\TagsOnRight

\def\Z{\Bbb Z}

\def\l{\left}
\def\r{\right}
\def\b{\bigg}
\def\bg{\bigg}
\def\({\b(}
\def\[{\b[}
\def\){\b)}
\def\]{\b]}

\def\t{\text}
\def\f{\frac}
\def\mo{\roman{mod}}
\def\em{\emptyset}
\def\se {\subseteq}
\def\sp {\supseteq}
\def\sm{\setminus}

\def\eq{\equiv}

\def\ls{\leqslant}
\def\gs{\geqslant}
\def\al{\alpha}

\def\la{\lambda}
\def\Proof{\noindent{\it Proof}}
\def\Remark{\noindent{\it Remark}}

\def\Ack{\noindent {\bf Acknowledgment}}
\hbox {Final version (Sept. 21, 2004) for J. Number Theory.}
\medskip
\title On the range of a covering function\endtitle
\author Zhi-Wei Sun\endauthor
\affil Department of Mathematics (and Institute of Mathematical Science)\\
Nanjing University, Nanjing 210093, The People's Republic of China\\
{\it E-mail}: {\tt zwsun$\@$nju.edu.cn}
\\Homepage: {\tt http://pweb.nju.edu.cn/zwsun}\endaffil

\abstract Let $\{a_s(\mo\ n_s)\}_{s=1}^k\ (k>1)$ be a finite system of
residue classes with the moduli $n_1,\ldots,n_k$ distinct.
By means of algebraic integers we show that the range of the covering function
$w(x)=|\{1\ls s\ls k:\, x\eq a_s\ (\mo\ n_s)\}|$ is not contained in any residue
class with modulus greater one. In particular, the values of
$w(x)$ cannot have the same parity.
\endabstract
\thanks 2000 {\it Mathematics Subject Classification}. Primary 11B25;
Secondary 11A07, 11A25, 11B75, 11R04.
\newline\indent
Supported by the Distinguished Youth Foundation
and the Key Program of NSF in P. R. China.
\endthanks
\endtopmatter
\document
\heading {1. Introduction}\endheading

  For $a\in\Z$ and $n\in\Z^+=\{1,2,3,\ldots\}$, let $a(n)$ stand for
the residue class $\{x\in\Z:\, x\eq a\ (\mo\ n)\}$. A finite
system
$$\{a_s(n_s)\}_{s=1}^k\ \ (k>1)\tag1.1$$
of residue classes is said to be a {\it cover} of $\Z$ if
$\bigcup_{s=1}^ka_s(n_s)=\Z$.

  The concept of cover of $\Z$ was introduced by P. Erd\H os ([E50]) in the
  early 1930s, who was particularly interested in those covers (1.1)
  with the moduli $n_1,\ldots,n_k$ distinct. By Example 3 of the author [S96],
  if $n>1$ is odd then
 $$\{1(2), 2(2^2),\ldots,2^{n-2}(2^{n-1}),2^{n-1}(n),2^{n-1}
 2(2n),\ldots,2^{n-1}n(2^{n-1}n)\}$$
  forms a cover of $\Z$ with distinct moduli.
  Covers of $\Z$ have been studied by various researchers (cf. [G04] and [PJ])
  and many surprising applications have been found
  (see, e.g. [F02], [S00], [S01] and [S03b]).

  Here are two major
  open problems concerning covers of $\Z$ (see sections E23, F13 and F14 of [G04]
  for references to these and other conjectures).

 \proclaim{Erd\H os--Selfridge Conjecture} Let $(1.1)$ be a cover
 of $\Z$ with distinct moduli. Then $n_1,\ldots,n_k$ cannot be all
 odd.\endproclaim

\proclaim{Schinzel's Conjecture} If $(1.1)$ is a cover of $\Z$,
then there is a modulus $n_t$ dividing another modulus $n_s$.
\endproclaim

For system (1.1), the function $w:\Z\to \Z$
given by
$$ w(x)=|\{1\ls s\ls k:\ x\in a_s(n_s)\}|\tag1.2$$
is called its {\it covering function}. Obviously $w(x)$
is periodic modulo the least common multiple $N=[n_1,\ldots,n_k]$
of the moduli $n_1,\ldots,n_k$.

Now we list some known results concerning the covering
function $w(x)$.

 (i) The arithmetic mean of $w(x)$ with $x$ in a period equals
 $\sum_{s=1}^k1/n_s$.

  (ii) (Z. W. Sun [S95, S96]) The covering function $w(x)$ takes its minimum on every set of
 $$\bg|\bg\{0\ls\theta<1:\, \sum_{s\in I}\f{m_s}{n_s}-\theta\in\Z
 \ \t{for some}\ I\se\{1,\ldots,k\}\bg\}\bg|$$
 consecutive integers, where $m_1,\ldots,m_k$ are given
 integers relatively prime to $n_1,\ldots,n_k$ respectively.

 (iii) (Z. W. Sun [S03a]) The maximal value of $w(x)$ can be written in the form
 $\sum_{s=1}^km_s/n_s$ with $m_1,\ldots,m_k\in\Z^+$.

 (iv) (\v S. Porubsk\'y [P75]) If $n_1,\ldots,n_k$ are distinct,
 then $[n_1,\ldots,n_k]$ is the smallest positive period
 of the function $w(x)$.

 (v) (Z. W. Sun [S03a]) If $n_0\in\Z^+$ is a period of the function
 $w(x)$, then for any $t=1,\ldots,k$ we have
 $$\bg\{\sum_{s\in I}\f1{n_s}:\, I\se\{1,\ldots,k\}\sm\{t\}\bg\}
 \sp\bg\{\f r{n_t}:\,r\in\Z\ \t{and}\ 0\ls r<\f{n_t}{(n_0,n_t)}\bg\},$$
 where $(n_0,n_t)$ denotes the greatest common divisor of $n_0$ and $n_t$.

 (vi) (Z. W. Sun [S04]) The function $w(x)$ is constant if $w(x)$ equals a constant
 for $|S|$ consecutive integers $x$ where
 $$S=\bg\{\f r{n_s}:\ r=0,\ldots,n_s-1;\ s=1,\ldots,k\bg\}.$$

  In this paper we study the range of a covering function (via algebraic integers)
  for the first time.
  Proofs of the theorems below will be given in the next section.

  \proclaim{Theorem 1.1} Suppose that the range of the covering function
  of $(1.1)$ is contained in a residue class with modulus $m$.
  Then, for any $t=1,\ldots,k$ with $mn_t\nmid [n_1,\ldots,n_k]$,
  we have $n_t\mid n_s$ for some $1\ls s\ls k$ with $s\not=t$.
  \endproclaim

  \proclaim{Corollary 1.1} If the covering function $w(x)$ of $(1.1)$ is constant, then
  for any $t=1,\ldots,k$ there is an $s\not=t$ such that $n_t\mid n_s$,
  and in particular $n_k=n_{k-1}$ provided that $n_1\ls\cdots\ls n_{k-1}\ls n_k$.
 \endproclaim
 \Proof. Suppose that $w(x)=c$ for all $x\in\Z$.
 Choose an integer $m>[n_1,\ldots,n_k]$. As $c(m)$ contains the range of $w(x)$,
 the desired result follows from Theorem 1.1. \qed

 \medskip
 \Remark\ 1.1. When $(1.1)$ is a disjoint cover of $\Z$,
 i.e., $w(x)=1$ for all $x\in\Z$, the first part
 of Corollary 1.1 was given by B. Nov\'ak and \v S. Zn\'am [NZ] and
 the second part was originally obtained by H. Davenport, L. Mirsky, D. Newman
 and R. Rad\'o independently. Corollary 1.1 appeared in Porubsk\'y [P75].

 \proclaim{Corollary 1.2} Suppose that those moduli in $(1.1)$ which are maximal with respect to divisibility
 are distinct. Then $w(\Z)=\{w(x):\,x\in\Z\}$ cannot be contained in a residue class
 other than $0(1)=\Z$, i.e., for any prime $p$ there is an $x\in\Z$ with
 $w(x)\not\eq w(0)\ (\mo\ p)$. In particular, those $w(x)$ with $x\in\Z$
 cannot have the same parity.
 \endproclaim
\Proof. Assume that $w(\Z)$ is contained in a residue class with modulus $m\in\Z^+$.
For each modulus $n_t$ maximal with respect to divisibility, there is no $s\not=t$
such that $n_t\mid n_s$, thus $mn_t$ divides $N=[n_1,\ldots,n_k]$ by Theorem 1.1.
Since $N$ is also the least common multiple of those moduli $n_t$ maximal with respect to divisibility,
we must have $mN\mid N$ and hence $m=1$. This ends the proof. \qed

\medskip
 \Remark\ 1.2. In contrast with the Erd\H os--Selfridge conjecture,
 Corollary 1.2 indicates that {\bf if (1.1) is a cover of $\Z$ with distinct moduli then
 not every integer is covered by (1.1) odd times}.
 \medskip

 Here is another related result.
\proclaim{Theorem 1.2} Let $A=\{a_s(n_s)\}_{s=1}^k$
and $B=\{b_t(m_t)\}_{t=1}^l$ both have
distinct moduli. Then $A$ and $B$ are identical
provided that $w_A(x)\eq w_B(x)\ (\mo\ m)$
for all $x\in\Z$, where $w_A$ and $w_B$ are covering functions of $A$ and $B$ respectively,
and $m$ is an integer not dividing $N=[n_1,\ldots,n_k,m_1,\ldots,m_l]$.
\endproclaim

\Remark\ 1.3. In 1975 Zn\'am [Z75] extended a uniqueness theorem of S. K. Stein [St] as follows:
Under the condition of Theorem 1.2, we have $A=B$ if $w_A=w_B$.
This follows from Theorem 1.2 by taking $m>N$.
\medskip

 Theorem 1.1 can be refined as follows.

 \proclaim{Theorem 1.3} Let $\la_1,\ldots,\la_k\in\Z$ be weights assigned to the $k$ residue classes
 in $(1.1)$ respectively.
 Suppose that $n_0\in\Z^+$ is the smallest positive period of
  $w(x)=\sum_{1\ls s\ls k,\, n_s\mid x-a_s}\la_s$
  modulo $m\in\Z$, and
that $d\in\Z^+$ does not divide $n_0$ but
  $I(d)=\{1\ls s\ls k:\, d\mid n_s\}\not=\em$.
  Then, either $m$ divides $[n_1,\ldots,n_k]\sum_{s\in I(d)}\la_s/n_s$, or we have
  $$|I(d)|\gs|\{a_s\ \mo\ d:\, s\in I(d)\}|\gs \min\Sb 0\ls s\ls k\\s\not\in I(d)\endSb \f d{(d,n_s)}\gs p(d)\tag1.3$$
  where $p(d)$ denotes the smallest prime divisor of $d$.
  \endproclaim

 \Remark\ 1.4. Theorem 1.3 in the case $m=0$ was first obtained by the author [S91] in 1991,
 an extension of this was given in [S04].
 \medskip

 Instead of (1.1) we can also consider a finite system
 of residue classes in $\Z^n$ (cf. [S04]) and deduce $n$-dimensional versions of Theorems 1.1--1.3.

 \heading{2. Proofs of Theorems 1.1--1.3}\endheading

 \noindent{\it Proof of Theorem 1.1}.
 Without any loss of generality we assume that $0\ls a_s<n_s$ for $s=1,\ldots,k$.
 Set $N=[n_1,\ldots,n_k]$. Then
 $$\align\sum_{r=0}^{N-1}w(r)z^r
 =&\sum_{r=0}^{N-1}\sum\Sb 1\ls s\ls k\\n_s\mid a_s-r\endSb z^r
=\sum_{s=1}^k\sum\Sb 0\ls r<N\\r\in a_s(n_s)\endSb z^r
\\=&\sum_{s=1}^kz^{a_s}\sum_{0\ls q<N/n_s}(z^{n_s})^q
\\=&\sum\Sb 1\ls s\ls k\\z^{n_s}=1\endSb\f{N}{n_s}z^{a_s}
+(1-z^N)\sum\Sb 1\ls s\ls k\\z^{n_s}\not=1\endSb\f{z^{a_s}}{1-z^{n_s}}.
 \endalign$$

 Suppose that $w(r)=a+mq_r$ for each $r\in\Z$ where $a,q_r\in\Z$.
 If $\al\not\in\Z$ but $\al N\in\Z$, then
 $$\sum_{r=0}^{N-1}w(r)e^{2\pi i\al r}=m\sum_{r=0}^{N-1}q_re^{2\pi i\al r}$$
 and also
 $$\sum_{r=0}^{N-1}w(r)e^{2\pi i\al r}=\sum^k\Sb s=1\\\al n_s\in\Z\endSb \f{N}{n_s}e^{2\pi i\al a_s},$$
 therefore we have the following congruence
 $$ \sum^k\Sb s=1\\\al n_s\in\Z\endSb \f{N}{n_s}e^{2\pi i\al a_s}\eq 0\ (\mo\ m)\tag2.1$$
 in the ring of all algebraic integers.

 If $1\ls t\ls k$ and $n_t\mid n_s$ for no $s\in\{1,\ldots,k\}\sm\{t\}$, then
 by applying (2.1) with $\al=1/n_t<1$ we obtain that
$$\f{N}{n_t}e^{2\pi i a_t/n_t}\eq 0\ (\mo\ m)$$
and hence $m$ divides $N/n_t$ in $\Z$.

The proof of Theorem 1.1 is now complete. \qed

\bigskip

\noindent{\it Proof of Theorem 1.2}. Without any loss of generality, we assume that
$n_1>\cdots>n_k$ and $m_1>\cdots>m_l$.
As $w_A(x)-w_B(x)\eq0\ (\mo\ m)$ for all $x\in\Z$, by modifying the proof of Theorem 1.1 slightly,
 we find that if $\al\not\in\Z$ but $\al N\in\Z$ then
 $$\sum^k\Sb s=1\\\al n_s\in\Z\endSb \f N{n_s}e^{2\pi i\al a_s}
 -\sum^l\Sb t=1\\\al m_t\in\Z\endSb \f N{m_t}e^{2\pi i\al b_t}\eq0\ \ (\mo\ m).\tag2.2$$
 In the case $d=\max\{m_1,n_1\}>1$,  by applying (2.2) with $\al=1/d$
 and the hypothesis $m\nmid N$, we get that
 $m_1=n_1$ and
 $$\f Nd\l(e^{2\pi i a_1/d}-e^{2\pi i b_1/d}\r)\eq0\ (\mo\ m).$$
 If $a_1\not\eq b_1\ (\mo\ d)$,  then $z=1-e^{2\pi i(b_1-a_1)/d}$
 is a zero of the monic polynomial $(-1)^{d-1}P(1-x)\in\Z[x]$ where $P(x)=(1-x^d)/(1-x)=1+x+\cdots+x^{d-1}$,
 hence $z$ divides the constant term $P(1)=d$ of $P(1-x)$
 in the ring of algebraic integers.
 As $m$ does not divide $N$, we must have $a_1\eq b_1\ (\mo\ d)$ and so $a_1(n_1)=b_1(m_1)$.
 Now that
$$|\{1<s\ls k:\,x\in a_s(n_s)\}|\eq|\{1<t\ls l:\,x\in b_t(m_t)\}|\ \ (\mo\ m),$$
we can continue the above procedure to obtain that
$$a_2(n_2)=b_2(m_2),\ \ldots,\ a_{\min\{k,l\}}(n_{\min\{k,l\}})=b_{\min\{k,l\}}(m_{\min\{k,l\}}).$$
If $k\not=l$, say $k>l$, then $m\Z$ contains the range of the covering function of $\{a_s(n_s)\}_{s=l+1}^k$
and this contradicts Theorem 1.1 since $m\nmid [n_{l+1},\ldots,n_k]$
and $n_{l+1}>\cdots>n_k$. So $A=B$ and we are done. \qed

\bigskip

\noindent{\it Proof of Theorem 1.3}.
Let $N=[n_1,\ldots,n_k]$. Clearly $(n_0,N)\in n_0\Z+N\Z$ is also a period of $w(x)$ mod $m$,
so $(n_0,N)=n_0$ and hence $n_0\mid N$.
Observe that
$$\sum^k\Sb s=1\\x\in a_s(n_s)\endSb\la_s-\sum^{n_0-1}\Sb r=0\\x\in r(n_0)\endSb
w(r)\eq0\ (\mo\ m)$$ for each $x\in\Z$.
As in the proof of Theorem 1.1, if $c\in\Z$ and $d\nmid c$ then
 $$\sum^k\Sb s=1\\(c/d)n_s\in\Z\endSb\la_s\f{N}{n_s}e^{2\pi i\f cda_s}
 -\sum^{n_0-1}\Sb r=0\\(c/d)n_0\in\Z\endSb w(r)\f{N}{n_0}e^{2\pi i\f cdr}
\eq0\ \ (\mo\ m). \tag2.3$$
For any $c\in\Z^+$ divisible by none of those $d/(d,n_s)$ with $0\ls s\ls k$ and $s\not\in I(d)$, we have
$$d\mid cn_s\iff\f{d}{(d,n_s)}\mid c\iff d\mid n_s\iff s\in I(d),$$
therefore (2.3) yields that
$$\sum_{s\in I(d)}\la_s\f {N}{n_s}e^{2\pi i\f cda_s}\eq0\ \ (\mo\ m).$$

 Let
 $$R=\{0\ls r<d:\, a_s\eq r\ (\mo\ d)\ \t{for some}\ s\in I(d)\}$$
 and suppose that $|R|<\min_{0\ls s\ls k,\, s\not\in I(d)}d/(d,n_s)$. By the above,
 $$u_n:=\sum_{r\in R}c_r\l(e^{2\pi i\f rd}\r)^n\eq0\ (\mo\ m)
\ \ \t{for every}\ n=1,\ldots,|R|,$$
where $c_r=N\sum_{s\in I(d),\,a_s\in r(d)}\la_s/n_s\in\Z$.
As $\{u_n\}_{n\gs0}$ is a linear recurrence of order $|R|$ with characteristic polynomial
$\prod_{r\in R}(x-e^{2\pi ir/d})$ whose coefficients are algebraic integers, we have
$u_n\eq0\ (\mo\ m)$ for every $n=|R|+1,|R|+2,\ldots$. In particular,
$\sum_{r\in R}c_r=u_d\eq0\ (\mo\ m)$,
i.e., $m$ divides $N\sum_{s\in I(d)}\la_s/n_s$.
We are done. \qed

\bigskip

\Ack. The author thanks the referee for his constructive comments.


\widestnumber\key{S03b}
\Refs

\ref\key E50\by P. Erd\H os\paper On integers of the form $2^k+p$
and some related problems\jour Summa Brasil.
Math.\vol2\yr1950\pages113--123\endref

\ref\key F02\by M. Filaseta\paper
Coverings of the integers associated with an irreducibility theorem
of A. Schinzel \jour in: Number Theory for the Millennium
(Urbana, IL, 2000), Vol. II, pp. 1-24, A K Peters, Natick, MA,
2002\endref

\ref\key G04\by R. K. Guy\book
Unsolved Problems in Number Theory\publ 3rd Edition, Springer, New York, 2004\endref

\ref\key P75\by \v S. Porubsk\'y\paper Covering systems and generating functions
\jour Acta Arith.\vol 26\yr1975\pages223--231\endref

\ref\key PS\by\v S. Porubsk\'y and J. Sch\"onheim \paper Covering
systems of Paul Erd\"os: past, present and future \jour in: Paul
Erd\"os and his Mathematics. I (edited by G. Hal\'asz, L.
Lov\'asz, M. Simonvits, V. T. S\'os), Bolyai Soc. Math. Studies
11, Budapest, 2002, pp. 581--627\endref

\ref\key NZ\by B. Nov\'ak and \v S. Zn\'am
\paper Disjoint covering systems\jour Amer. Math. Monthly
\vol 81\yr 1974\pages 42--45\endref

\ref\key St\by S. K. Stein\paper Unions of arithmetic sequences
\jour Math. Ann.\vol 134\yr 1958\pages 289--294\endref

\ref\key S91\by Z. W. Sun\paper An improvement to the Zn\'am--Newman
result\jour Chinese Quart. J. Math.\vol 6\yr 1991\pages no. 3, 90--96\endref

\ref\key S95\by Z. W. Sun\paper Covering the integers by
arithmetic sequences
 \jour Acta Arith.\vol 72\yr1995\pages109--129\endref

\ref\key S96\by Z. W. Sun\paper Covering the integers by
arithmetic sequences {\rm II}
 \jour Trans. Amer. Math. Soc.\vol348\yr1996\pages4279--4320\endref

\ref\key S00\by Z. W. Sun\paper On integers not of the form $\pm
p^a\pm q^b$ \jour Proc. Amer. Math. Soc.\vol 128\yr 2000\pages
997--1002\endref

\ref\key S01\by Z. W. Sun\paper Algebraic approaches to periodic
arithmetical maps\jour J. Algebra\vol 240\yr
2001\pages723--743\endref

\ref\key S03a\by Z. W. Sun\paper On the function
$w(x)=|\{1\ls s\ls k:\, x\eq a_s\ (\mo\ n_s)\}|$
\jour Combinatorica\vol 23\yr 2003\pages 681--691\endref

\ref\key S03b\by Z. W. Sun\paper Unification of zero-sum problems,
subset sums and covers of $\Z$ \jour Electron. Res. Annnounc.
Amer. Math. Soc. \vol 9\yr 2003\pages 51--60\endref

\ref\key S04\by Z. W. Sun\paper Arithmetic properties of periodic
maps \jour Math. Res. Lett.\vol 11\yr 2004\pages 187--196\endref

\ref\key Z75\by \v S. Zn\'am\paper On properties of systems of arithmetic sequences
\jour Acta Arith.\vol 26\yr 1975\pages 279--283\endref

\endRefs
\enddocument